\documentclass[12pt]{article}
\usepackage{epsf,colordvi}
\newcommand{\ft}[2]{{\textstyle\frac{#1}{#2}}}
\newsavebox{\uuunit}
\sbox{\uuunit}
    {\setlength{\unitlength}{0.825em}
     \begin{picture}(0.6,0.7)
        \thinlines
        \put(0,0){\line(1,0){0.5}}
        \put(0.15,0){\line(0,1){0.7}}
        \put(0.35,0){\line(0,1){0.8}}
       \multiput(0.3,0.8)(-0.04,-0.02){12}{\rule{0.5pt}{0.5pt}}
     \end {picture}}
\newcommand {\unity}{\mathord{\!\usebox{\uuunit}}}
\newcommand  {\Rbar} {{\mbox{\rm$\mbox{I}\!\mbox{R}$}}}

\newcommand {\Cbar}
    {\mathord{\setlength{\unitlength}{1em}
     \begin{picture}(0.6,0.7)(-0.1,0)
        \put(-0.1,0){\rm C}
        \thicklines
        \put(0.2,0.05){\line(0,1){0.55}}
     \end {picture}}}
\newcommand{\Ka}{K\"ahler}

\def\Im{{\rm Im ~}}
\def\Re{{\rm Re ~}}

\newcommand{\rmi}{{\rm i}}

\newcommand{\sinprod}[2]{\mbox{$\langle #1 , #2 \rangle$}}
\newcommand{\Poin}{Poincar\'e}
\newcommand{\Sp}[1]{\mbox{$Sp\left( #1,\Rbar \right) $}}

\csname @addtoreset\endcsname{equation}{section}

\begin{document}

\begin{titlepage}
\begin{flushright}
KUL-TF-2000/04\\
math.DG/0002122
\end{flushright}
\vspace{.5cm}
\begin{center}
\baselineskip=16pt
{\LARGE    Special K\"{a}hler geometry  
}\\
\vfill
{\large Antoine Van Proeyen $^{\dagger}$, 
  } \\
\vfill
{\small
 Instituut voor Theoretische Fysica, Katholieke
 Universiteit Leuven,\\
Celestijnenlaan 200D B-3001 Leuven, Belgium
}
\end{center}
\vfill
\begin{center}
{\bf Abstract}
\end{center}
{\small
The geometry that is defined by the scalars in couplings of
Einstein--Maxwell theories in $N=2$ supergravity in 4 dimensions is
denoted as special \Ka\ geometry. There are several equivalent definitions,
the most elegant ones involve the symplectic duality group. The
original construction used conformal symmetry, which immediately
clarifies the symplectic structure and provides a way
to make connections to quaternionic geometry and Sasakian manifolds.
}\vspace{2mm} \vfill
Contribution to the Proceedings of the Second Meeting on Quaternionic
Structures in Mathematics and Physics, Rome 6-10 September 1999.\vfill
\hrule width 3.cm
{\footnotesize
\noindent $^\dagger$ Onderzoeksdirecteur, FWO, Belgium }
\end{titlepage}
\section{Introduction}
In the previous workshop in this series on quaternionic geometry, B.
de Wit and me gave talks \cite{prtrquat} on the classification of quaternionic
homogeneous spaces \cite{HomogSpecial}. Results in special geometry had lead to new
homogeneous quaternionic spaces. We have discussed on this topic
further with D. Alekseevsky and V. Cort\'{e}s\footnote{V. Cort\'{e}s made our results
more accessible to the mathematical audience \cite{Cortes}.} and realised that it would
be useful to have a definition of special \Ka\ geometry that does not
refer to the constructions of supersymmetric actions. The text of the
proceedings was a first step in that direction.
Meanwhile, in 1994, the second superstring revolution took place. The
main issue was that theories which were previously thought as
different, are recognized as perturbations around `vacua' of a master
theory. Essential for that are the duality relations which make the
connections between the different descriptions. The first example was
provided by Seiberg and Witten \cite{SeiWit}. They used a model with
$N=2$ supersymmetry in 4 dimensions with vector multiplets, being
multiplets involving Maxwell fields. Special \Ka\ geometry \cite{dWVP} is defined
by the couplings of the scalars in the locally supersymmetric theory, i.e.\
in the coupled Einstein--Maxwell theory. The model used by Seiberg--Witten
thus involves a similar geometry, which has been called rigid special
\Ka\ geometry \cite{STG}, as it appears in rigid supersymmetry. The structure of
that geometry was important for the obtained results. In particular
the analyticity properties of fields in these theories allowed them
to find exact solutions.

\begin{figure}[htb]
\begin{center}
\epsfbox[0 350 288 527]{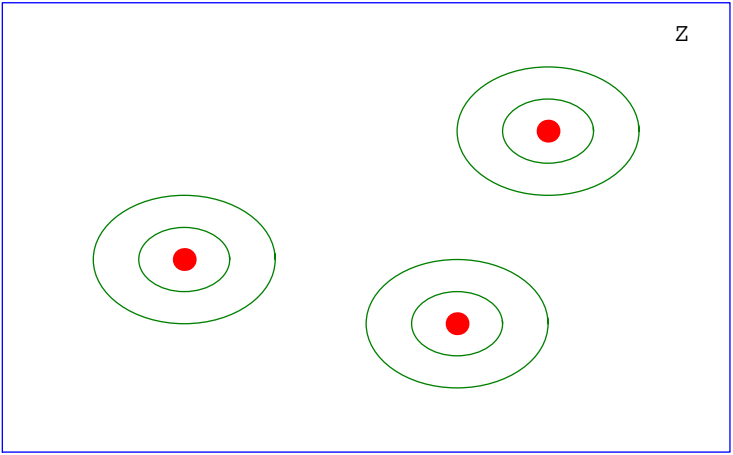}
\caption{\it The moduli space with 3 singularities.}\label{fig:holz}
\end{center}
\end{figure}

The so-called vector multiplets in $d=4$, $N=2$ supersymmetry are
multiplets with spins $(0,0,\ft12,\ft12,1)$, the latter being the
vector providing the Maxwell theory. The scalars are moduli, whose values
parametrize the different vacua. The two (real) scalars in a
multiplet can be combined to a complex one, and the supersymmetry
will indeed provide a complex structure. As will become clear below,
the structure of special \Ka\ geometry implies holomorphicity of the
resulting field equations. Then the result of Seiberg--Witten is based
on the fact that singularities and the asymptotic
behaviour determine exact answers. The singularities, see
figure~\ref{fig:holz},
are points around which a classical limit can be considered. The
theory allows perturbation expansions around these points. Each one
leads classically to a different theory, but there is only one full
quantum theory. The singular points form a family of inequivalent
vacua.

These developments motivated us to look for a definition of special
geometry independent of supersymmetry. A first step in that direction
had meanwhile be taken by Strominger \cite{strom}. He had in mind the
moduli spaces of Calabi--Yau spaces. His definition is already based
on the symplectic structure, which we also have emphasized. However,
being already in the context of Calabi--Yau moduli spaces, his
definition of special \Ka\ geometry omitted some ingredients that are
automatically present in any Calabi--Yau moduli space, but have to be
included as necessary ingredients in a generic definition. Another important step
was obtained in \cite{CDFVP}. Before, special geometry was connected
to the existence of a holomorphic prepotential function $F(z)$. The
special \Ka\ manifolds were recognized as those for which the \Ka\
potential can be determined by this prepotential, in a way to be
described below. However, in \cite{CDFVP} it was found that one can
have $N=2$ supergravity models coupled to Maxwell multiplets such
that there is no such prepotential. These models were constructed
by applying a symplectic transformation to a model with prepotential.
This fact raised new questions: are all the models without
prepotential symplectic dual to models with a prepotential? Can one
still define special \Ka\ geometry starting from the models with a
prepotential? Is there a more convenient definition which does not
involve this prepotential? These questions have been answered in
\cite{Craps:1997gp}, and are reviewed here.

Section~\ref{ss:ingredients} introduces some ingredients. I give some
elements of the algebraic context of $N=2$ supersymmetry, and how the
geometric quantities are encoded in the action. Then I show the
emergence of symplectic transformations in the actions with vector
fields coupled to scalars. Rigid $N=2$ supersymmetry and the
associated rigid special \Ka\ geometry is discussed in
section~\ref{ss:rigid}. Section~\ref{ss:N2sg} will then discuss the
supergravity case. For that, it is useful to look first at the
conformal group, as a formulation from that perspective will show more
structure, in particular it clarifies the role of the symplectic
transformations, and gives the connection with Sasakian manifolds.
This is the central section where the definitions, their equivalence
and some examples are discussed. The special \Ka\ manifolds appear in
moduli spaces of Riemann surfaces for the rigid version and in those
of Calabi--Yau manifolds for the local version. That is illustrated
in section~\ref{ss:realRSCY}. A summary is
given in section~\ref{ss:summary}. We briefly discuss there also the
usage of the same construction methods for quaternionic geometry as
recently applied in \cite{QuatConf}.

\section{Ingredients} \label{ss:ingredients}
For supersymmetry in 4-dimensional spacetime, the fermionic charges
should belong to a spinor representation of $SO(3,1)$. Therefore,
in the minimal supersymmetric case, the supercharges have 4 real components.
This minimal situation is called $N=1$. Field theory allows
realizations up to $N=8$ supersymmetry, i.e.\ with 32 real
supercharges. Special \Ka\ geometry appears in the context of
$N=2$ supersymmetry. The 8 real spinor supercharges are
denoted as $Q_\alpha ^i$, where $\alpha =1,\ldots ,4$ and $i=1,2$.
They satisfy the anticommutation rule
\begin{equation}
\left\{ \Red{Q^i_\alpha} ,\Red{Q^j_\beta} \right\} =
\gamma^\mu_{\alpha \beta } \RoyalPurple{P_\mu}\delta ^{ij}\,,
\label{anticQ}
\end{equation}
thus involving the translation operator $P_\mu $ in 4-dimensional
spacetime. There are representations with spins
\begin{eqnarray}
(0,0,0,0,\ft12,\ft12):  & \mbox{hypermultiplet} & \mbox{quaternionic scalars} \nonumber\\
(0,0,\ft12,\ft12,1): & \mbox{vector multiplet}   & \mbox{complex scalars} \nonumber\\
(1,\ft32,\ft32,2):& \mbox{supergravity} &\,,
\label{reprN2}
\end{eqnarray}
where I have indicated their names and the types of scalars. The
quaternionic and complex structures are guaranteed by the
supersymmetry.

The ingredients of the geometry are found in the action. In general,
having scalars $\Blue{\phi ^i(x)}$, vectors with field strength
$\Blue{{\cal F}_{\mu\nu}^I(x)}$, and possibly a non-trivial spacetime metric
$\Plum{ g_{\mu \nu }(x)}$, the bosonic kinetic part of the action has
the general form
\begin{eqnarray}
 S&=& \int d^4x \Plum{\sqrt{g}g^{\mu \nu }} \Blue{\partial _\mu \phi ^i
\partial _\nu \phi ^j} \Red{G_{ij}(\phi )}\nonumber\\
&& +\ft14\Plum{\sqrt{g}g^{\mu \rho }g^{ \nu \sigma }}
 \Red{(\Im {\cal N}_{IJ})(\phi )}\Blue{{\cal F}_{\mu\nu}^I
{\cal F}_{\rho \sigma }^ J}
-\ft \rmi 8 \Red{(\Re {\cal N}_{IJ})(\phi )}
\Plum{\varepsilon^{\mu\nu\rho\sigma}}\Blue{{\cal F}_{\mu\nu}^I
{\cal F}_{\rho\sigma}^J}\nonumber\\
&& + \ldots\,. \label{bosAct}
\end{eqnarray}
$\Red{G_{ij}(\phi )}$ is identified as the metric of the manifold of
scalars. The complex symmetric matrix
$\Red{{\cal N}_{IJ}}$ determines the kinetic terms of the vectors,
and its meaning will be clarified below.

Supersymmetry relates bosons and fermions, e.g.\ for the scalars
\begin{equation}
   \delta \Blue{\phi ^i(x)} = \overline{\Red{\epsilon }}\ \Blue{\chi^i
   (x)}\,,
\label{delSUSY}
\end{equation}
where $\Red{\epsilon }$ are the supersymmetry parameters and $\Blue{\chi^i
(x)}$ are the fermions. In the context of local supersymmetry the parameters
depend on spacetime, and we thus have
\begin{equation}
  \delta \Blue{\phi ^i(x)} = \overline{\Red{\epsilon (x)}}\,\Blue{\chi^i(x)}\,.
\label{delSUGRA}
\end{equation}
In order to have an action invariant under these local symmetries,
one needs connection fields, which are the gravitini for the supersymmetry.
Due to the algebra (\ref{anticQ}) this should be related to local
translations, i.e.\ general coordinate transformations, whose
connection field is the (spin 2) graviton.

A prerequisite to understand the following development, is the understanding
of the meaning of the \Red{symplectic transformations}. These are
the duality symmetries of 4 dimensions, the generalizations of the
Maxwell dualities. They were first discussed in
\cite{dual}. Consider the kinetic terms of the vector fields as in
(\ref{bosAct}) with $I=1, ..., m$. $\Red{{\cal N}_{IJ}}$ are coupling
constants or functions of scalars. One defines (anti)selfdual
combinations as
\begin{equation}
  \Blue{{\cal F}^\pm _{\mu\nu}}= \ft12 \left( \Blue{{\cal F}_{\mu\nu}}  \pm \ft12
\varepsilon_{\mu\nu\rho\sigma}\Blue{{\cal F}^{\rho\sigma}}\right)\,.
\label{selfdual}
\end{equation}
The conventions\footnote{The Levi--Civita symbol has $\varepsilon _{0123}=\rmi$.}
are such that the complex conjugate of $\Blue{{\cal F}^+}$
is $\Blue{{\cal F}^-}$. Defining
\begin{equation}
  \Blue{G}_{+I }^{\mu\nu}\equiv 2\rmi\frac{\partial{\cal L}}
  {\partial {\cal F}^{+I }_{\mu\nu}}=
\Red{{\cal N}_{I J }}\Blue{{\cal F}^{+J \,\mu\nu}}\,,
\label{defG}
\end{equation}
the Bianchi identities and field equations can be written as
\begin{eqnarray}
\partial^\mu \Im \Blue{{\cal F}^{+I }_{\mu\nu} }&=&0\ \ \ \ \ {\rm Bianchi\
identities}\nonumber\\
\partial_\mu \Im \Blue{G_{+I }^{\mu\nu}} &=&0\ \ \ \ \  {\rm Equations\  of\
motion.}\label{BianchiField}
\end{eqnarray}
This set of equations is invariant under \RawSienna{$GL(2m,\Rbar)$}:
\begin{equation}
\Blue{\pmatrix{\widetilde{\cal F}^+\cr \widetilde G_+\cr}}={\cal S}
\Blue{\pmatrix{{\cal F}^+\cr G_+\cr}} =
\pmatrix{A&B\cr C&D\cr}  \Blue{ \pmatrix{{\cal F}^+\cr G_+\cr}}\,.
\end{equation}
In order that this transformation be consistent with (\ref{defG}), we
should have
\begin{eqnarray}
&&\Blue{\widetilde G^+}=(C+D\Red{{\cal N}})\Blue{{\cal F}^+}=
(C + D\Red{{\cal N}})(A+B\Red{{\cal N}})^{-1} \Blue{\widetilde {\cal F}^+} \nonumber\\
&&\rightarrow \mbox{
\parbox[t]{8.5cm}
{\fbox{$\Red{\widetilde{\cal N}} = (C + D\Red{{\cal N}})(A+B\Red{{\cal
N}})^{-1}$}}}\,.\label{sympltrN}
\end{eqnarray}
However, this matrix should remain symmetric,
 $\Red{\widetilde {\cal N}}=\Red{\widetilde {\cal N}^T}$, which
 implies that
\begin{equation} {\cal S}=\pmatrix{A&B\cr C&D\cr} \in \RawSienna{Sp(2m,\Rbar)}\,,
\end{equation}
as the explicit condition is
\begin{equation}  {\cal S}^T  \Omega   {\cal S}   =  \Omega  \qquad\mbox{where}\qquad
\Omega=\pmatrix{0&\unity \cr -\unity &0\cr}\,.
\label{Scond}
\end{equation}
Thus the remaining transformations are real symplectic ones in
dimension $2m$, where $m$ is the number of vector fields.

In the following we will denote by symplectic vectors, those vectors
$\RawSienna{V}$ such that its symplectic transformed is
 $\RawSienna{\widetilde V} = {\cal
S}\RawSienna{V}$. The prime example is thus
$\RawSienna{V}= \Blue{ \pmatrix{{\cal F}^+\cr G_+\cr}}$.
An invariant inner product of symplectic vectors is defined by
\begin{equation}
  \langle \RawSienna{V} , \RawSienna{W} \rangle \equiv
 \RawSienna{V^T} \Omega \RawSienna{W}\,.
\label{innerProdSymplV}
\end{equation}

The important properties for the matrix $\Red{{\cal N}}$ is that it
should be symmetric and $\Im\Red{{\cal N}}<0$ in order to have
positive kinetic terms. These properties are
\Red{preserved under symplectic transformations} defined by (\ref{sympltrN}).
\section{Rigid special \Ka\ geometry} \label{ss:rigid}
As mentioned in the introduction, the `rigid' special \Ka\ geometry
is the geometric structure encountered in rigid $N=2$ supersymmetry
in 4 dimensions. This supersymmetry has as field representations
multiplets with spins $(0,0,0,0,\ft12,\ft12)$, the hypermultiplet,
and multiplets with spin $(0,0,\ft12,\ft12,1)$, the vector multiplet.
For the former, the scalar field geometry is based on quaternions,
and is a hyper-\Ka\ structure. Here, we will consider the vector
multiplets, for which the scalars combine to complex fields, whose
geometry is \Ka ian. A natural description for such multiplets uses
$N=2$ superspace, that is an extension of usual spacetime (with
points labelled by $x$) by fermionic coordinates $\theta $, such that
the superspace is a representation of the superalgebra. The vector
multiplets are then described by superfields $\Green{\Phi^A (x,\theta )}$
that satisfy some constraints, restricting the way in which they
depend on the $\theta $. The result is some superfield
\begin{equation}
  \Green{\Phi ^A(x,\theta )}= \Blue{X^A(x)}+ \bar\theta \Blue{\chi _\alpha ^A(x)} +\bar \theta \gamma ^{\mu \nu }
   \theta  \Blue{{\cal F}_{\mu \nu }(x)} + \ldots\,,
\end{equation}
where the lowest components $\Blue{X^A}$ are complex fields. $A=1,\ldots
,n$ labels different vector multiplets. To build an action, one
integrates a general holomorphic function $F$ over one half of the $\theta
$ variables (the chiral superspace). The above mentioned constraints
have, between other restrictions, restricted the superfields to
depend only on this chiral superspace. With
\begin{equation}
S=\int d^4x\int d^4\theta\ F(\Green{\Phi})\ + c.c.\,,
\end{equation}
one obtains that the scalars have a metric of \Ka ian type:
\begin{eqnarray}
\Red{G _{A \bar B}(X,\bar X)} & = &\partial_A\partial_{\bar B}\Magenta{K(X,\bar X)}  \nonumber\\
\Magenta{K(X,\bar X)}&=&i (\bar F_A(\bar X) X^A-F_A(X) \bar X^A)\nonumber\\
\Red{{\cal N}_{AB}}&=& F_{AB}\,,
\end{eqnarray}
where the latter defines the kinetic term of the vectors as in
(\ref{bosAct}). Further,
$F_A(X) =\frac{\partial}{\partial X^A}
F(X)$ or  $\bar F_A(\bar X)=\frac{\partial} {\partial \bar X^A}\bar F(\bar
X)$, $F_{AB}(X) =\frac{\partial}{\partial X^A}\frac{\partial}{\partial X^B}
F(X)$.

The equations of motion turn out to be those equations that determine that
$F_A(\Phi)$  satisfy the same superfield constraints as $\Phi^A$.
Comparing with (\ref{BianchiField}), the superfield constraints on $\Phi
^A$ contain the first equations (Bianchi identities) while the same
superfield equations on $F_A$ contain the second line.

It is therefore appropriate to combine the superfield in a
`symplectic vector'
\begin{equation}
\RawSienna{\tilde\Phi}=\pmatrix{\Phi^A\cr F_A(\Phi)}\qquad
\begin{array}{l}
\mbox{chiral superfields which} \\ \mbox{satisfy extra constraints.}
\end{array}
\end{equation}
The scalars, i.e.\ the $\theta =0$ part of this vector form also a
symplectic vector
\begin{equation}  \RawSienna{V}\equiv   \pmatrix{X^A\cr F_A(X)}
 \qquad\mbox{is a symplectic vector.}
\end{equation}

A further improvement is to allow general coordinates. So far, we
parametrize the scalars as $X^A$, which are special coordinates (occurring in
the superfields). We can, however, allow arbitrary coordinates \cite{CdAF}
$z^\alpha$ with  $\alpha=1,\ldots ,n$. Then the special coordinates
are holomorphic functions of the $z^\alpha $, i.e.\ $X^A (z^\alpha )$, such that
$e_\alpha^A\equiv \partial_\alpha X^A(z)$ is invertible.

Now we have all the ingredients to give definitions \cite{Craps:1997gp}.

\Blue{{\bf Definition 1 of rigid special \Ka\ geometry.}}

A rigid special \Ka\ manifold is an $n$-dimensional \Ka\ manifold with on any chart $n$ holomorphic
functions $X^A(z)$ and a \Plum{holomorphic function $F(X)$} such that
\begin{equation} \Magenta{K(z,\bar z)}=\rmi \left( X^A\frac{\partial}{\partial\bar X^A}\bar F(\bar X)-  \bar X^A
\frac{\partial}{\partial X^A} F(X)\right)\,. \end{equation}
On overlap of charts these functions should be related by
(inhomogeneous) symplectic transformations $ISp(2n,\Rbar)$:
\begin{equation}
\left( \begin{array}{c}
 X \\ \partial F \end{array}\right)_{(i)} = e^{\rmi c_{ij}} M_{ij}
\left( \begin{array}{c}
 X \\ \partial F\end{array}\right)_{(j)}+ b_{ij}\,, \label{ISpn}
\end{equation}
with
\begin{equation} c_{ij}\in \Rbar\ ;\qquad M_{ij} \in \Sp{2n}\ ;\qquad b_{ij} \in
\Cbar^{2n} \,,
\end{equation}
satisfying the cocycle condition on overlaps of 3 charts.

There is, however, a second definition of rigid special \Ka\
manifolds, which is based on the symplectic structure, rather than on
the prepotential.

\Blue{{\bf Definition 2 of  rigid special \Ka\ geometry.}}

A \Ka\ manifold is the base manifold of a $U(1)\times I\Sp{2n}$  bundle.
A holomorphic section $\RawSienna{V}(z)$ defines the \Ka\ potential
by
\begin{equation}
\Magenta{K(z,\bar z)}=\rmi \sinprod {\RawSienna{V}}{\RawSienna{{\bar V}}}\,,
\end{equation}
and it should satisfy the constraint
\begin{equation}
 \sinprod {\partial_\alpha \RawSienna{V}}{\partial_\beta \RawSienna{V}} =0\,.
\label{extrarigidconstraint}
\end{equation}

One can show that the prepotential exists locally, but it is thus not
essential for the definition. In rigid special geometry the choice of
definition is rather a question of esthetics. However, in the local
case, it will be important to have the analogue of the second
definition available. The kinetic matrix for the vectors is
\begin{equation}
  \Red{{\cal N}_{AB}}=\left( \partial_\alpha F_A(z)\right)
  e_B^\alpha\,.
\label{Ndef2rigid}
\end{equation}
The condition (\ref{extrarigidconstraint}) guarantees that this
matrix is symmetric.
Finally, let us remark that the symplectic metric $\Omega $ should in general
not assume the canonical form (\ref{Scond}), but can be an arbitrary
non-degenerate real antisymmetric matrix. However, in order to
distinguish $X^A$ and $F_A$ components, and thus to write a
prepotential, one should bring it to this canonical form.

\section{$N=2$ supergravity and special \Ka\ geometry}\label{ss:N2sg}

In this section we introduce the `local' special \Ka\ geometry, which
is the one generally denoted as special \Ka\ geometry. It is this one
which was found in \cite{dWVP}, and has most interesting
applications. It was introduced in the context of supergravity. To
explain its structure, it is useful to consider again its origin.

To describe a supergravity theory, there are several methods. One of
them is the introduction of a \Blue{{\bf superspace}}. This formalism
shows a \Plum{lot of structure} of the theory. It is very transparent
for rigid supersymmetry. However, in its local version, necessary for
supergravity, there appear a lot of \Plum{extra superfield symmetries}.
These symmetries are an artifact of the formalism. They
have to be gauge-fixed to obtain the physical theory.

\Blue{{\bf Superconformal tensor calculus}} is in-between. Also here
\Plum{extra gauge symmetries} occur, and these are in fact the symmetries
of the superconformal group, the basic ingredient of the formalism.
The experience tells us that
these symmetries are the relevant ones to display the structure of the
theory, but this formalism does not have the many other symmetries present
in the superspace approach. It turns out that we just remain with those that
are useful to get insight in complicated formulae. Also for the calculation
of the action, the superconformal symmetries are just appropriate to simplify
the construction.
This is \Red{particularly interesting} in our case.  The superconformal
tensor calculus gives the proper setup for the \Red{symplectic}
(duality-adapted) formulation.

The \MidnightBlue{idea} is to start by constructing an action invariant under
superconformal group.  Then, one
choose gauges for the extra gauge invariances of the superconformal group, such
that the remaining theory has just the super-\Poin\ symmetries.

The formalism can be used for theories in various dimensions and amount of
supersymmetry. Let us review here the structure for 4 dimensions with
8 real supersymmetry generators ($N=2$). The superconformal group
contains first of all the \Maroon{conformal group}
(translations, Lorentz rotations, dilatations and
special conformal transformations). This group is
 \Maroon{$S0(4,2)=SU(2,2)$}.
The supersymmetries should sit in a spinor representation of this group.
This singles out the supergroup $SU(2,2|2)$, which means essentially
that the group can be represented by matrices of the form
\begin{equation}
\pmatrix{\Maroon{SU(2,2)} & \Red{SUSY} \cr \Red{SUSY} & \Green{SU(2)\times
U(1)}}\,.\label{supergroup}
\end{equation}
The off-diagonal blocks are the fermionic symmetries. The diagonal
blocks are the bosonic ones. They split up in the above-mentioned
conformal group and an `\Green{$R$}-symmetry group', $\Green{SU(2)\times
U(1)}$. This extra group plays an important role:
\begin{itemize}
  \item  the gauge connection of $U(1)$ will be the \Red{\Ka\ curvature}.
  It acts on the manifold of scalars in vector multiplets,
  \item the gauge connection of $SU(2)$ promotes the hyper\Ka\ manifold of
  hypermultiplets to a \Red{quaternionic} manifold.
\end{itemize}

As we neglect here the hypermultiplets, we have to consider the basic
supergravity multiplet and the vector multiplets. The physical
content that one should have (from representation theory of the
super-\Poin\ group) can be represented as follows:
\begin{equation}
\begin{array}{ccccccccc}
\multicolumn{3}{c}{SUGRA}&\qquad&&\multicolumn{3}{c}{vector m.}&\\
&2&                      &      &&&&&\\
\ft32&&\ft32             &      &&&&&\\
    &1&                  &      &&&1&&\Red{\rightarrow n+1}\\
    & &                  &      &+ \Green{n}\ *\quad&\ft12 & & \ft12 &\\
    & &                  &      &           &   \Green{0}  & &  \Green{0}&
\end{array}
\end{equation}
The supergravity sector contains the graviton, 2 gravitini and a
so-called graviphoton. That spin-1 field gets, by coupling to $n$
vector multiplets, part of a set of \Red{$n+1$} vectors, which will be
uniformly described by the special \Ka\ geometry. The scalars appear
as \Green{$n$} complex ones \Green{$z^\alpha$, with $\alpha=1,\ldots , n$}.

To describe this, we start with \Red{$n+1$ superconformal vector multiplets} with scalars
\Blue{$X^I$ with $I=0,\ldots ,n$}. The action is determined by a
holomorphic function $F(X)$. Compared with the rigid case, there is
one additional requirement. The conformal invariance requires $F(X)$ to be
homogeneous of weight~2, where the $X$ fields carry weight~1. These
scalar fields transform also under a local $U(1)$ symmetry.

The obtained metric is a cone \cite{coneGary,QuatConf}. To see this, one
splits the $n+1$ complex variables $\{\Blue{X}\}$ in
$\{\rho ,\theta , \Green{z^\alpha} \}$
\begin{itemize}
  \item
$r$ is scale which is a gauge degree of freedom for translations
  \item
$\theta $ is the $U(1)$ degree of freedom;
  \item the $n$ complex variables $\Green{z^\alpha }$.
\end{itemize}
The metric now takes the form
\begin{equation}
  ds^2=dr^2 +\ft1{18}r^2 \left[A+ d\theta + \rmi\left(\partial _\alpha K(z,\bar z)\, dz^\alpha -
  \partial _{\bar \alpha }K(z,\bar z)\,d\bar z^{\bar \alpha }\right)\right] ^2+
  r^2\partial _\alpha \partial _{\bar \alpha }K(z,\bar z)\, dz^\alpha d\bar z^{\bar \alpha
  }\,,
\label{dsSasakian}
\end{equation}
where $A$ is the one-form gauging the $U(1)$ group, and $K(z,\bar z)$
is a function of the holomorphic prepotential $F(X)$, to be explained
below.
With $A=0$, this defines the cone over a Sasakian manifold. However,
in supergravity, the field equation of $A$ implies that it is a
composite field, given by (minus) the other parts of the second term of
(\ref{dsSasakian}). With fixed $\rho $ (gauge fixing the superfluous
dilatations), the remaining manifold is \Ka, with the \Ka\ potential
determined by $F(X)$. That gives the special \Ka\ metric.

Let us explain this now in more detail, using at the same time more
of the symplectic formalism. The dilatational gauge fixing (the
fixing of $r$ above), is done by the condition
\begin{equation}
  X^I\bar F_I(\bar X)-\bar  X^I F_I(X)=\rmi\,.
\label{dilatgauge}
\end{equation}
This condition is chosen in order to decouple kinetic terms of the
graviton from those of the scalars. Using again symplectic vectors
\begin{equation}
  \RawSienna{V}=\pmatrix{X^I\cr F_I}\,,
\label{Vlocal}
\end{equation}
this can be written as the condition on the symplectic inner product:
\begin{equation}
  <\RawSienna{V},\RawSienna{\bar V}>=\rmi\,.
\label{VbarVisi}
\end{equation}
To solve this condition, we define
\begin{equation}
  \RawSienna{V}=e^{\Magenta{K(z,\bar z)}/2}\RawSienna{v}(z)\,,
\label{Vtov}
\end{equation}
where $v(z)$ is a holomorphic symplectic vector,
\begin{equation}
\RawSienna{v}(z)= \pmatrix{Z^I(z)\cr \frac{\partial }{\partial Z^I}F(Z)
}\,.
\label{vz}
\end{equation}
The upper components here are arbitrary functions (up to conditions
for non-degeneracy), reflecting the freedom of choice of coordinates
$z^\alpha $. The \Ka\ potential is
\begin{equation}
  e^{-\Magenta{K(z,\bar z)}}=-\rmi \sinprod {\RawSienna{v}}{\RawSienna{{\bar
  v}}}\,.
\label{Kahlerv}
\end{equation}
The kinetic matrix for the vectors is given by
\begin{equation} \Red{{\cal N}_{IJ}}= \pmatrix{F_I&{\cal D}_{\bar \alpha} \bar F_I(\bar  X)}
\pmatrix{X^J&{\cal D}_{\bar \alpha} \bar X^J}^{-1}\,,
\label{Nsugra}
\end{equation}
where the matrices are $(n+1)\times (n+1)$ and
\begin{equation}
  {\cal D}_{\bar \alpha} \bar F_I(\bar  X)=\partial_{\bar \alpha} \bar F_I(\bar  X)+\ft12
  (\partial_{\bar   \alpha}\Magenta{K})\bar F_I(\bar  X)\,,\qquad
{\cal D}_{\bar \alpha} \bar X^J=\partial_{\bar \alpha} \bar X^J+\ft12
  (\partial_{\bar   \alpha}\Magenta{K})\bar X^J\,.
\label{Dbarbar}
\end{equation}

Before continuing with general statements, it is time for an \Red{example}.
Consider the prepotential $F=-\rmi X^0X^1$. This is a model with $n=1$.
There is thus just one coordinate $z$. One has to choose a
parametrization to be used in the upper part of (\ref{vz}). Let us
take a simple choice: $Z^0=1$ and $Z^1=z$. The full symplectic vector
is then (as e.g.\ $F_0(Z)=-\rmi Z^1(z)$)
\begin{equation}
  \RawSienna{v}=\pmatrix{Z^0\cr Z^1\cr F_0\cr F_1}=\pmatrix{1\cr z\cr -\rmi z\cr
  -\rmi}\,.
\label{vexample}
\end{equation}
The \Ka\ potential is then directly obtained from (\ref{Kahlerv}),
determining the metric:
\begin{equation}
e^{-\Magenta{K}}=2(z+\bar z)\ ; \qquad \Red{g_{z{\bar z}}}
=\partial_z\partial_{{\bar z}}\Magenta{K}=
(z+{\bar z})^{-2}\,.\label{KaMetricExample}
\end{equation}
The kinetic matrix for the vectors is diagonal. From (\ref{Nsugra})
follows
\begin{equation}
\Red{{\cal N}}=\pmatrix{-\rmi z&0\cr 0& -\rmi\ft1z\cr}\ .
\end{equation}
Therefore the action contains
\begin{equation} e^{-1}{\cal L}_1= -\ft12 \Re\left[z \left( \Blue{F_{\mu\nu}^{+0}}\right) ^2
+ z^{-1} \left( \Blue{F_{\mu\nu}^{+1}}\right) ^2 \right]\,.
\label{firstL1ex}
\end{equation}
The domain of positivity for both metrics is $\Re z>0$.

We formulate again two definitions, the first using the prepotential,
and the second one using only the symplectic vectors.

\Blue{{\bf Definition 1 of (local) special \Ka\ geometry.}}

A special \Ka\ manifold is
an \Plum{$n$-dimensional Hodge-\Ka\ manifold} with on any chart
\Plum{$n+1$ holomorphic functions $Z^I(z)$} and a
\Magenta{holomorphic function $F(Z)$, homogeneous of second degree}, such
that, with (\ref{vz}), the \Ka\ potential is given by
\begin{equation}
e^{-K(z,\bar z)}=-\rmi \sinprod v{\bar v}\,,
\end{equation}
and on overlap of charts, the $v(z)$ are connected by
\Plum{symplectic transformations} $Sp(2(n+1),\Rbar)$
and/or  \Plum{\Ka\ transformations}.
\begin{equation} v(z) \rightarrow e^{f(z)}{\cal S} v(z)\,.
\end{equation}

\Blue{{\bf Definition 2 of  (local) special \Ka\ geometry.}}

A special \Ka\ manifold is an $n$-dimensional \Ka--Hodge manifold,
that is the base manifold of a  $Sp(2(n+1))\times U(1)$ bundle. There
should exist a \Red{holomorphic section $v(z)$} such that the \Ka\
potential can be written as
\begin{equation}
e^{-K(z,\bar z)}=-\rmi \sinprod v{\bar v}\,,
\end{equation}
and it should satisfy the condition
\begin{equation}
\sinprod {{\cal D}_\alpha v}{{\cal D}_\beta v} =0\,.
\label{DavDbv0}
\end{equation}

Note that the latter condition guarantees the symmetry of ${\cal
N}_{IJ}$. This condition did not appear in \cite{strom}, where the
author had in mind Calabi--Yau manifolds. As we will see below,
in those applications, this condition is automatically fulfilled.
For $n>1$ the condition can be replaced by the equivalent condition
\begin{equation}
\sinprod {{\cal D}_\alpha v}{v} =0\,.
\label{Davv0}
\end{equation}
For $n=1$, the condition (\ref{DavDbv0}) is empty, while
(\ref{Davv0}) is not. In \cite{sympl} it has been shown that models
with $n=1$ not satisfying (\ref{Davv0}) can be formulated.

The appearance of `Hodge' manifold in the definitions refers to a
global requirement. The $U(1)$ curvature should be of even integer
cohomology. This has been considered first in \cite{BaggerWitten},
and for an explanation on the normalization, one can consult
\cite{Craps:1997gp}. Note that in the mathematics literature `Hodge'
refers to integer cohomology. Here, however, the presence of fermions
makes the condition stronger by a factor of two: one needs even
integers.

Let us come back to the example, on which we will perform a
symplectic mapping:
\begin{equation}
\RawSienna{\widetilde v}={\cal S}\RawSienna{v}=
\pmatrix{1&0&0&0\cr 0&0&0&-1\cr 0&0&1&0\cr 0&1&0&0\cr }
\RawSienna{v}= \pmatrix{1\cr \rmi\cr -\rmi z\cr z}\,.
\label{tildevExample}
\end{equation}
After this mapping, $z$ is not any more a good coordinate for
$(\widetilde Z^0,\widetilde Z^1)$, the upper two components of the
symplectic vector $z$. This means that the symplectic vector can not
be obtained from a prepotential. We can not obtain the symplectic
vector from a form (\ref{vz}). No function
$\widetilde F(\widetilde Z^0,\widetilde Z^1)$ exists.
Therefore, the first definition is not applicable. However,
nothing prevents us from using the second definition. The \Ka\ metric
is still the same, (\ref{KaMetricExample}), and one can again compute
the vector kinetic matrix, either directly from (\ref{Nsugra}), as
the denominator is still invertible, or from (\ref{sympltrN}):
\begin{equation}
\Red{\widetilde{\cal N}}= (C+D\Red{ {\cal N}})(A+B \Red{{\cal N}})^{-1}=
-\rmi X^1(X^0)^{-1}\unity=-\rmi z \unity \,.
\end{equation}
In this parametrization, the action is thus
\begin{equation} e^{-1}{\cal L}_1= -\ft12 \Re \left[z \left(\Blue{ F_{\mu\nu}^{+0}}\right) ^2
+z \left(\Blue{ F_{\mu\nu}^{+1}}\right) ^2 \right]\,.
\end{equation}
This action is not the same as the one before, but is a `dual
formulation' of the same theory, being obtained from (\ref{firstL1ex})
by a duality transformation. The straightforward construction in
superspace or superconformal tensor calculus does not allow to
construct actions without a superpotential. However, in \cite{sympl}
it has been shown that the field equations of these models can also
be obtained from the superconformal tensor calculus. One just has to
give up the concept of a superconformal invariant action.

It is thus legitimate to ask about the \Blue{equivalence of the two
definitions}. Indeed, we saw that in some cases definition 2 is satisfied,
but one can not obtain a prepotential $F$.
However, that example, as others in \cite{CDFVP}, was obtained from performing a
symplectic transformation from a formulation where the prepotential
does exist. In \cite{Craps:1997gp} it was shown that this is true in
general. If definition 2 is applicable, then there exists a
symplectic transformation to a  \Plum{basis such that $F(Z)$ exists}.
Note, however, that in the way physical problems are handled, the
existence of formulations without prepotentials is important. Going
to a dual formulation, one obtains a formulation with different symmetries
in perturbation theory. The example that we used here appears in a
reduction to $N=2$ of two versions of $N=4$ supergravity, known respectively
as the `$SO(4)$ formulation' \cite{N4SO4}
and the `$SU(4)$ formulation' of pure $N=4$ supergravity \cite{CSF}.

Finally let us note that we still could apply (\ref{Nsugra}) because
the matrix
\begin{equation}\pmatrix{X^I&{\cal D}_\alpha \bar X^I} \end{equation}
is \Red{always invertible} if the metric $g_{\alpha \bar \alpha }=
\partial _\alpha \partial _{\bar \alpha }K(z,\bar z)$ is positive
definite.
Therefore, the inverse exists, and ${\cal N}_{IJ}$ can be constructed.
However, the matrix
\begin{equation}\pmatrix{X^I&{\cal D}_\alpha  X^I} \end{equation}
is not invertible in the formulation (\ref{tildevExample}).
\Blue{\textbf{If}} that matrix
is \Blue{invertible}, \Blue{then} a \Blue{prepotential} exists
\cite{Craps:1997gp}.
\section{Realizations in moduli spaces of Riemann surfaces and Calabi--Yau manifolds}
\label{ss:realRSCY}
The realizations of special \Ka\ geometry that are mostly studied in
physics these days, are the moduli spaces of Riemann surfaces for
the rigid case, and those of Calabi--Yau 3-folds for the local case.

First, consider the Hodge diamond of Riemann surfaces, listing the
number of non-trivial (anti)holomorphic $(p,q)$ forms:
\begin{center}
\begin{tabular}{ccc}
& $h^{00}=1$      &  \\
\Plum{$h^{10}=g$ & & $h^{01}=g$}\\
 &$h^{11}=1$&
\end{tabular}
\end{center}
\Blue{Rigid special \Ka\ }geometry is obtained for the moduli spaces
of such \Blue{Riemann surfaces} when we consider
\begin{itemize}
\item
with \Red{$n$ complex moduli} $z^\alpha$
\item $n\leq g$ \Red{ holomorphic 1-forms}
 $\gamma_\alpha$  ($\alpha=1,\dots\, n$)
\item
\Red{$2n\ $  cycles $c_\Lambda$} that form a complete basis for
1-cycles for which $\int_c \gamma_\alpha\neq 0$.
\end{itemize}
In this situation
\begin{equation}
\gamma_\alpha (z)   =\partial_\alpha \lambda(z) + d\eta_\alpha(z)\,,
\end{equation}
where $\lambda(z)$ is a meromorphic 1-form with zero residues.
The symplectic formulation of rigid special \Ka\ geometry is obtained
with as symplectic vector the vector of periods of $\lambda $ over
the chosen cycles:
\begin{equation} \RawSienna{V=\int_{c_\Lambda} \lambda}\,.  \end{equation}
The intersection matrix of the cycles plays the role of the
symplectic metric. This type of realizations was used in
Seiberg--Witten models. The general features have been discussed in
\cite{Craps:1997gp}.

To obtain local special \Ka\ manifolds, one considers the
\Green{moduli space of Calabi--Yau 3-folds}.
In this case the Hodge diamond of the manifold is
\begin{center}
\begin{tabular}{ccccccc}
 & &  & $h^{00}=1$      &   &  & \\
 & &0 &        & 0 &  & \\
 &0&  &$h^{11}=m$&   & 0& \\
\Plum{$h^{30}=1$& &$h^{21}=n$&  &$h^{12}=n$&&$h^{03}=1$}\\
 &0&  &$h^{22}=m$&   & 0& \\
 & &0 &        & 0 &  & \\
 & &  & $h^{33}=1$      &   &  &
\end{tabular}
\end{center}
These manifolds have $h^{21}=n$ complex structure moduli, which play
the role of the variables \Green{$z^\alpha$} of the previous section.
There are
\Plum{$2(n+1)$ 3-cycles $c_\Lambda$}, with intersection matrix
$Q_{\Lambda\Sigma}=c_\Lambda\cap  c_\Sigma$. The canonical form is
obtained with so-called $A$ and $B$ cycles, and then $Q$ takes the form
of $\Omega $ in (\ref{Scond}). Symplectic vectors are identified
again as vectors of integrals over the $2(n+1)$ 3-cycles:
\begin{equation}
\RawSienna{v=\int_{c_\Lambda} \Omega^{(3,0)}}\,,\qquad
{\cal D}_\alpha v= \int_{c_\Lambda} \Omega_\alpha ^{(2,1)}\,.
\end{equation}
$\Omega^{(3,0)}$ is the unique $(3,0)$ form that characterizes the Calabi--Yau
manifold. $\Omega_\alpha ^{(2,1)}$ is a basis of the $(2,1)$ forms, determined
by the choice of basis for $z^\alpha $.
That these moduli spaces give rise to special \Ka\ geometry became
clear in \cite{cy2skg}. Details on the relation between the geometric
quantities and the fundamentals of special \Ka\ geometry have been
discussed in \cite{fresoriabook,Craps:1997gp}.

The defining equations of special \Ka\
geometry are automatically satisfied. E.g.\ one can easily see
how the crucial equation (\ref{DavDbv0}) is realized:
\begin{eqnarray}
\left<{\cal D}_\alpha v,{\cal D}_\beta v \right>&=& \int_{c_\Lambda}
\Omega^{(2,1)}_{(\alpha)}
\cdot Q^{\Lambda\Sigma}\cdot  \int_{c_\Sigma}\Omega^{(2,1)}_{(\beta)}
\nonumber\\ &=&
\int_{CY}\Omega^{(2,1)}_{(\alpha)}  \wedge \Omega^{(2,1)}_{(\beta)}= 0\,.
\end{eqnarray}
The \Maroon{symplectic transformations} correspond now to changes of
the basis of the cycles used to construct the symplectic vectors. The
statement that a formulation with a prepotential can always be
obtained in special \Ka\ geometry by using a symplectic
transformation, can now be translated to the statement that the
geometry can be obtained from a \Maroon{prepotential for some choice of cycles}.

Finally, it is interesting that singularities of Calabi--Yau
manifolds may be used to obtain a `rigid limit'. Indeed, in \cite{KLocalRigid}
it is shown how Calabi--Yau manifolds that are $K3$ fibrations can
be reduced near the singularity to fibrations of ALE manifolds. Then
the special geometry of the moduli space of the Calabi--Yau manifold
reduces to the rigid special geometry with \Ka\ potential determined
by the ALE manifold. This mechanism is considered further in
\cite{cyk3}. There it has been shown how the \Ka\ potential of
special geometry approaches the one of rigid special geometry, and
how the periods of the local theory behave around the singular points
and thus around the rigid limit. In the superstring theory this
allows to get the gravity corrections to the rigid theory, which can
be used for applications \cite{CacZanDenef}.

\section{Summary and connection with quaternionic manifolds}\label{ss:summary}
\Red{Special \Ka\ }geometry is defined by the couplings of
$N=2$ supersymmetric theories (`rigid' special \Ka) or supergravity theories
((local) special \Ka) with vector multiplets. There are several ways to
describe the geometry. We discussed two ways:
\begin{itemize}
\item by using a prepotential function
\item by symplectic vectors and constraints
\end{itemize}
In \Blue{{\it rigid}} special \Ka\ geometry, these are completely
equivalent. In the \Blue{{\it local}} theory, all special \Ka\ manifolds
can be obtained from a prepotential, but in some cases that involves
a duality transformation. Therefore not all actions can be described by the
prepotential.

\Blue{{\it Rigid}} special \Ka\ geometry is realised by moduli spaces of
certain Riemann manifold. That construction is not straightforward,
and involves a choice of cohomology subspace and moduli.
The \Blue{{\it local}} special \Ka\ geometry appears in the moduli space
of Calabi--Yau threefolds. In this case the construction is
straightforward. For a particular Calabi--Yau manifold one includes
all the moduli. In this way a
clear geometrical interpretation of the building blocks of special
geometry is obtained. \Maroon{Duality transformations} correspond then to
a change of the basis of cycles. A prepotential does exist at least
for a suitable choice of basis of the cycles.

Note, however, that not all special \Ka\ manifolds can be obtained
as realizations in moduli spaces. E.g.\ the homogeneous manifolds,
treated in \cite{prtrquat,HomogSpecial} are not obtained in this way.

In the First Meeting on Quaternionic
Structures in Mathematics and Physics, 5 years ago, we have shown
\cite{prtrquat,HomogSpecial} how
\Red{homogeneous} special \Ka\ spaces are related by the \Blue{$c$-map} to
\Blue{homogeneous quaternionic spaces} and by the \Green{$r$-map} to
\Green{homogeneous `very special' real spaces}.
The construction of special \Ka\ geometry that we have outlined here
can be used as well for the quaternionic spaces (and for the real ones). In a recent work
\cite{QuatConf} it has been shown how the conformal tensor calculus
can be applied to obtain the actions based on the quaternionic
spaces (actions for `hypermultiplets').
The scalars are the lowest components of superfields (or superconformal
multiplets) \Blue{$A_i^\alpha $}, with $i=1,2$ and $\alpha =1,\ldots ,
2(r+1)$ with a reality condition. The $A_i^\alpha $ can be considered
as \Blue{$Sp(1)\times Sp(r+1)$ sections}. Again the number of multiplets
($r+1$) is one more than the number of physical multiplets ($r$) that
we will obtain. We thus start with
$4(r+1)$ scalars. One of those will be a scale degree of freedom\footnote{When
vector multiplets and hypermultiplets are simultaneously coupled, there is one
overall dilatational gauge degree of freedom. An auxiliary field of
the superconformal gauge multiplet gives a second relation, such that
as well the compensating field of the vector multiplet as that of the
hypermultiplet are fixed.}, three are $SU(2)$ degrees of freedom,
the second part of the $R$-symmetry as was mentioned after
(\ref{supergroup}), and the remaining ones form \Red{$r$ quaternions}.
As in the metric of the vector multiplets, there is a connection to
Sasakian manifolds. Putting the gauge fields of the $SU(2)$ invariance to
zero, rather than using their field equations, one obtains a
 \Green{3-Sasakian manifold}. This is related to the talk of Galicki
 in the meeting 5 years ago \cite{Galicki}.

\medskip
\section*{Acknowledgments.}

\noindent
I am grateful to S. Vandoren for a useful discussion on the relation with
Sasakian manifolds. Most of this
review treats work done in collaboration with B. de Wit, B. Craps, F.
Roose and W. Troost, and I learned a lot from the
discussions with them. This work was
supported by the European Commission TMR programme ERBFMRX-CT96-0045.

\end{document}